# Anxiety profiles and protective factors: A latent profile analysis in children

Irene C. Mammarella[1], Enrica Donolato[1], Sara Caviola[2], & David Giofrè[3]

[1] *Department of Developmental Psychology and Socialization, University of Padova, Italy*

[2] *Department of Psychology, University of Cambridge, UK*

[3] *Department of Natural Sciences and Psychology, Liverpool John Moores University, UK*

Correspondence to:

Irene C. Mammarella

Department of Developmental Psychology and Socialization

University of Padova

Via Venezia 8

35131 Padova, Italy

e-mmail: irene.mammarella@unipd.it

**Research highlights**

- Different anxiety profiles and the influence of personal protective factors (self-concept and resilience) in schoolchildren were examined

- Three different profiles (low, moderate and high risk) were identified by a latent profile analysis

- Protective factors, such as self-concept and resilience, were differently related to anxiety

- Self-concept was lower the higher the anxiety risk profile

- Resilience only decreased in association with the high anxiety risk profile

**Anxiety profiles and protective factors: A latent profile analysis in children**

**Abstract**

The current study investigated the presence of different anxiety profiles in schoolchildren in order to understand whether Mathematics and Test Anxiety are a manifestation of a general form of anxiety, or the expression of specific forms of anxiety. Moreover, we also examined the influence of personal protective factors. The results of a latent profile analysis, conducted on 664 children attending grades 3 to 6, clearly identified three different profiles distinguished on the basis of the level of general, test and mathematics anxiety. Protective factors, such as self-concept and resilience, were differently related to anxiety: The former was clearly lower when the risk profile was higher, whereas students were able to maintain a certain level of resilience up to an average risk of developing forms of anxiety. The implications of these findings may lead to the development of specific intervention programs aimed at reducing students' anxiety and fostering self-concept and resilience.

**Anxiety profiles and protective factors: A latent profile analysis in children**

Anxiety is an aversive motivational state that occurs in situations where the level of perceived threat to the individual is high (Eysenck & Calvo, 1992). Different forms of anxiety have been described in the literature but one of the most prominent is general anxiety (GA): this refers to an individual's disposition to worry about many different events, behaviors or personal abilities of everyday life, together with a difficulty in controlling these worries (Eysenck & Calvo, 1992). While GA refers to a general condition of anxiety, other more specific forms of anxiety have also been described and have received a great deal of attention in the literature. Specifically, math anxiety, commonly defined as a feeling of tension, apprehension, or fear which may interfere with one's performance of mathematical tasks (Richardson & Suinn, 1972), and test anxiety (TA), a psychological, physical, or behavioral reaction to worry cognitions regarding potential failure in achievement/school assessment situations (Zeidner, 1998).

In the extant literature, however, it is not clear whether MA and TA are a manifestation of a general form of anxiety, or the expression of specific forms of anxiety. As a consequence, mathematics anxiety (MA) and test anxiety (TA) have often been investigated separately. Their effects have been examined in relation either to academic performance (Hill et al., 2016; Putwain, Daly, Chamberlain, & Sadreddini, 2015; Roick & Ringeisen, 2017) or to other forms of anxiety, such as general anxiety (Carey, Devine, Hill, & Szűcs, 2017). To the best of our knowledge, different forms of anxiety and any factors protecting against it have not been investigated together, as part of the same study. Protective factors are conditions or attributes, such as strengths, resources, supports, that operate in different domains of functioning (i.e., individual, school, family) and help the individuals to foster competence, promote successful development, and mitigate all the conditions associated with a higher likelihood of negative outcomes, or risk factors (Dekovic, 1999).

In the present study, we considered personal protective factors (i.e., self-concept and resilience), and distinguished between their general and academic effects. To the best of our knowledge, this is the first contribution on the association of different forms of anxiety with personal protective factors in school-age children.

### The relationship between general-, mathematics-, or test anxiety with personal protective factors

The results of previous meta-analyses suggest that GA is more closely related to TA than to MA (Hembree, 1988, 1990; Ma, 1999). Although the relationship between GA and MA is not particularly strong, GA seems to have a systematic effect on MA and on mathematics achievement. For example, Hill and collegues (2016) found that partialling out the effect of GA reduced the significant negative relationship between MA and mathematics achievement in primary- and middle-school students. Regarding the relationship between GA with either self-concept or resilience in college students, previous studies revealed an association between negative affect responses to aversive situations and lower levels of self-esteem or self-concept ( e.g., Moreland & Sweeney, 1984; Smith & Petty, 1995; see also Lowe, Papanastasiou, Deruyck, & Reynolds, 2005; Lowe, Peyton, & Reynolds, 2007). Notably, in graduate students, the correlation between GA and self-concept was moderate and negative. Benetti and Kambouropoulos (2006), using path analyses, examined the influence of resilience and GA on self-concept in undergraduate students and investigated the mediating role of positive and negative effects. Their findings indicated that positive and negative affects significantly mediated the influence of resilience and GA on self-concept, respectively with any significant direct effects between GA, resilience and self-concept.

Coming to MA, an extensive body of literature documents that not only cognitive factors, but also low self-confidence in math or negative attitudes to math teachers, are related to poor performance in math classes (Ashcraft, Kirk & Hopko, 1998; Ma, 1999; Maloney & Beilock, 2012; Mammarella, Caviola, Giofrè, & Borella, 2017; Mammarella, Hill, Devine, Caviola, & Szűcs, 2015). As for the role of personal protective factors, several studies have analyzed the relationship between academic self-concept and MA. In a cross-cultural study examining 15 years-olds from 41 countries, Lee (2009) found that MA and math self-concept are inversely related to one another. In a recent study, Justicia-Galiano, Martín-Puga, Linares, and Pelegrina (2017) investigated whether working memory and math self-concept mediate the relationship between MA and math performance among school-age children. Their results indicated that both working memory and math self-concept, as mediators, contributed to explaining the relationship between MA and mathematics achievement. Intriguingly, other researchers have suggested that MA is antecedent to self-concept and self-esteem (Ahmed, Minnaert, Kuyper & van der Werf, 2012). This would mean that MA can promote negative academic self-concepts regarding math abilities (Ashcraft & Kirk, 2001; Wu, Barth, Amin, Malcarne, & Menon, 2012).

Tests can trigger another type of academic anxiety, called test anxiety. Classic measures of TA have used "worry" and "emotionality" to distinguish between its cognitive and affective-physiological aspects, respectively (Morris, Davis, & Hutchings, 1981). More recent measures have considered other features of TA too, such as specific autonomic bodily reactions of anxiety, social evaluation and cognitive interference (Benson, Moulin-Julian, Schwarzer, Seipp, & El-Zahhar, 1992; Lowe, Grumbein, & Raad, 2011; Wren & Benson, 2004). Several studies have attested to the relationship between TA and protective factors like academic self-concept and self-esteem, which describe the individual's self-perceived ability in academic situations, with a strong impact on a student's TA (Goetz, Preckel, Zeidner, &

Schleyer, 2008; Hembree, 1988). Research has identified poor self-directed competence beliefs as a major cause of TA (Bong & Skaalvik, 2002; Putwain & Daniels, 2010). Bandalos, Yates, and Thorndike-Christ (1995) also found that academic self-concept was negatively related to TA, and mediated the link between prior experience and anxiety in undergraduate students. The relationship between resilience and TA has been examined too. According to Martin and Marsh (2006), resilient pupils can perform better under pressure, in anxiety-inducing testing situations, by maintaining a strong belief in their ability (which supports their effort and persistence), or by controlling worries likely to interfere with their performance. Putwain, Nicholson, Connors, and Woods (2013) examined whether TA mediated the relations between resilience and performance in a high-stake tests, after controlling for prior abilities in students attending grade 6. They found a significant indirect relationship between resilience and test performance, that was mediated by TA. A particular form of resilience, called academic buoyancy, also seems to be related to TA. Academic buoyancy is a relatively new construct emerging from the resilience literature, and defined as the student's response to academic challenges and pressures (Martin & Marsh, 2006, 2008, 2009). Academic buoyancy is associated with positive academic outcomes, including lower levels of TA and higher levels of persistence and confidence (Martin, Colmar, Davey, & Marsh, 2010; Putwain, Connors, Symes, & Douglas-Osborn, 2012). According to the self-referral model of TA devised by Zeidner and Matthews (2005), academic buoyancy could be seen as a factor protecting against maladaptive forms of coping and avoidance. For instance, buoyant students may have experienced failure in the past, or may expect to fail in the future, but this does not give rise to attributions likely to reduce their effort or motivation simply to protect themselves. Putwain and colleagues (2012) also found that academic buoyancy explained a significant portion of the variance in all components of TA (ranging from $R^2=.13$ to $R^2=.23$), and was inversely related thereto.

**The present Study**

Earlier research tended to study the effects of GA, TA, or MA in isolation, with scarce consideration for the effects of personal protective factors such as self-concept and resilience. Carey and colleagues (2017) recently assessed different forms of anxiety by conducting a latent profile analysis on students in grade 4, or in grades 7 and 8. They identified four profiles in grade 4, ranging from low to high anxiety. This four-group solution also emerged on students in grades 7 and 8, but the profiles appeared more specific in this case, and were described as: low anxiety; general anxiety; academic anxiety (i.e., MA and TA); and high anxiety. Studying such latent profiles, it is interesting to see how distinct but related forms of anxiety appear within a population. Unlike simple correlations, latent profiles help us to identify heterogeneous subgroups that express certain anxiety patterns. For these reasons, in the present study we conducted a latent profiles analysis on a large group of $3^{rd}$- to $6^{th}$-grade primary-school children, considering their scores on GA, TA and MA. Our first aim was to test whether specific latent profiles of anxiety emerged between $3^{rd}$ - and $6^{th}$ –graders in order to better understand whether MA and TA emerged as unique forms of a more general apprehension (GA) or specific expressions of different form of anxiety, at least for these particular age groups. We expected to find different profiles of anxiety derived by the combination of GA, TA and MA, in agreement with the results of meta-analytical studies which showed moderate correlations among these variables (Hembree, 1988, 1990; Ma, 1999). We focused on a wide age range to shed more light on whether it is already possible to distinguish between different forms of anxiety (mainly general and academic) even in younger children. However, given that previous studies shown that TA and MA peak around grades 9 to 10, and do not change thereafter (Ashcraft & Moore, 2009; Hembree, 1990), we did not expect strong age-differences in our sample.

Another aim of our study was to test the role of personal protective factors such as self-concept and resilience in relation to latent profiles of different forms of anxiety. As mentioned previously, academic self-concept and resilience, or academic buoyancy, both seem to be inversely related to TA. While previous studies mainly investigated the relationship between academic self-concept or academic buoyancy and anxiety (and TA in particular), here we distinguished between general and academic personal protective factors. This is because we assumed that, just as it seems important to distinguish between general and academic forms of anxiety, so too could a distinction between general and academic personal protective factors help to clarify their potential relationship between different latent profiles of anxiety. We, thus, expected that general and academic personal protective factors would have a different influence on latent profiles of anxiety, in agreement with previous studies showing a negative relation between TA and academic buoyancy (Martin et al., 2010).

## Method

### Participants

The study originally involved 666 children. However, some of the data from two children were missing and we decided to exclude these data from the analyses. Our final sample included 664 schoolchildren ($M_{age}$ = 9.20 years, SD = 1.13 years), 52.6% males and 47.4% females, attending primary school, in grades 3 to 6 (with 184, 206, 166, and 108 children in grades, 3, 4, 5 and 6, respectively). Seven different schools were involved and recruited through well-established contacts of the first Author. The children came from middle-class families[1], and were attending schools in urban areas of north-east Italy: 5% of the sample was composed of African ethnicity and 3% of other ethnicities. The study was approved by the Ethics Committee on Psychology Research at the University of Padova,

Italy. Parental consent was obtained. Children with intellectual disabilities or neurological/genetic disorders were not included in the study.

## Materials

### General anxiety

The *Revised Children's Manifest Anxiety Scale: Second Edition* (RCMAS-2; Reynolds & Richmond, 2012) is a self-report questionnaire for detecting general anxiety in children from 6 to 19 years old. It comprises 49 items with a yes/no response format. The questionnaire consists of two scales, the Total Anxiety scale and the Defensiveness scale, and three anxiety subscales, Social Anxiety, Physiological Anxiety, and Worry.. Items assessed, for example, whether the child often presents stomach ache, whether s/he is worried that his/her classmates could make fun of him/her and whether the child feels nervous when things don't go as s/he wants. The total anxiety score (Cronbach α = .89) was used for the present study. As reported in the manual (see Reynolds & Richmond, 2012), the RCMAS-2 presents good convergent validity (r =.61) and test-retest reliability (r= .71) for the total anxiety score.

### Academic anxiety

*Math-anxiety.* The *Abbreviated Math Anxiety Scale* (AMAS; Hopko, Mahadevan, Bare, & Hunt, 2003) is a self-report tool developed to assess MA. The Italian version of the AMAS for children was used in the present study (Caviola, Primi, Chiesi, & Mammarella, 2017). Like the original English version, the questionnaire contains 9 items that describe different situations involving math. Children were asked to assess each item in terms of how anxious they would be, scoring them on a 5-point scale ranging from 1 "strongly agree" to 5 "strongly disagree", the higher the number the greater their mathematical anxiety. A good internal consistency (Cronbach's α = .77) was observed. Finally, strong convergent validity (*r*= .85) and excellent 2-weeks rest-retest reliability (*r*= .85) were found in adults (see Hopko et al., 2003).

*Test-anxiety*. The *Test Anxiety Questionnaire for Children* (TAQ-C; Donolato, Marci, Altoè, & Mammarella, submitted) a new self-report for assessing TA in children was used. The questionnaire was translated and readapted from the Children's Test Anxiety Scale (CTAS; Wren & Benson, 2004) from which three subscales were used: thoughts (e.g., "I think I'm going to get a bad grade"); off-Task Behaviors (e.g., "I play with my pencil"); and autonomic reactions (e.g., "My heart beats fast"). The social derogation subscale from the FRIEDBEN Test Anxiety Scale for Adolescents (FTA; Friedman & Bendas-Jacob, 1997) was considered too (e.g., "I am worried that all my friends will get high scores in the test and only I will get low ones"). The self-report consisted of 37 items. Participants gave their answers on a 4-point Likert scale (from "almost never" to "almost always"). In the current sample, Cronbach's α was .84. Moreover, good convergent validity (*r* ranging from .71 to .77) and test-retest validity (*r*=.74) for the total score were found (Donolato et al., submitted).

**General protective factors**

*Self-Concept - Competence scale* (SC-Competence scale; Italian translation, Bracken, 2003). The *Competence* subscale is part of the Multidimensional Self-Concept Scale (MSC) that is developed for children between 9 and 19 years old. Specifically, the *Competence* subscale is composed by 25 items used to measure general protective factors (e.g., "I trust in myself"). Participants responded on a 4-point Likert scale from "absolutely true" to "absolutely false". The total score was calculated as recommended in the manual (Bracken, 2003), so higher scores corresponded to a greater degree of competence self-concept (Cronbach's α = .79). As reported in the manual (see Bracken, 2003), the subscale is characterized by good 4-weeks test-retest reliability (*r*= .76) as well as good convergent validity for the total MSC score (*r*=.73).

The *Ego-Resiliency scale* (ER; Block & Kremen, 1996) is brief inventory used to detect a set of traits relating to general resourcefulness, strength of character, and flexibility of functioning in response to different environmental demands, which reflects an individual's capacity for appropriate self-regulation. The scale has been used both with adults and children (see Caprara, Steca & De Leo, 2003). Participants used a 4-point scale from 1 ("does not apply at all") to 4 ("applies very strongly") to score 14 items, including "I quickly get over and recover from being startled", and "I enjoy dealing with new and unusual situations". Good internal consistency (Cronbach's $\alpha$ = .75) and 2-years test-retest reliability ($r$= .39) scores were observed in children and adolescents (see Caprara, Steca & De Leo, 2003).

**Academic protective factors**

*Self-Concept – Academic scale* (SC-Academic scale; Italian translation, Bracken, 2003). The *academic* subscale is part of the Multidimensional Self-Concept Scale (MSC) that is devised for children between 9 and 19 years old. The subscale was used to assess academic protective factors (e.g., "I am good in mathematics") and contains 25 items. Participants gave their answers on a 4-point Likert scale, from ("absolutely false") to ("absolutely true"). The total score was calculated as explained in (Bracken, 2003), so higher scores corresponded to higher levels of academic self-concept (Cronbach's $\alpha$ = .87). The *academic* subscale showed good test-retest reliability ($r$= .81) as well as good convergent validity regarding the total MSC score ($r$=.73), as reported by the authors (see Bracken, 2003).

*Academic Buoyancy Scale* (ABS; Martin & Marsh, 2008a, 2008b) is a brief inventory for assessing the ability to face and deal successfully with academic difficulties and challenges typical of ordinary school life (e.g., poor grades, competing deadlines, exam pressure, difficult schoolwork). The questionnaire is composed of 4 items (e.g. "I'm good at dealing with setbacks at school") rated from 1 ("strongly disagree") to 7 ("strongly agree"), which have been translated for the present study in Italian. The scale showed good internal consistency (Cronbach's α was .82) and test-retest reliability ($r = .67$) as reported by the authors (see Martin & Marsh, 2008a).

**Procedure**

The children were tested in a collective session in their classroom lasting approximately 60 minutes, during which the self-report questionnaires were administered in a fixed pseudorandomized order (SC-Competence scale, SC-Academic scale, RCMAS-2, AMAS, ER, TA, ABS). Questionnaires were presented one at the time and students had a short break after the completion of each questionnaire. Questionnaires were administrated by the second Author, a second-year PhD student in psychology. Children were not compensated for their participation, which is a standard practice in Italy.

**Statistical analyses**

R was used in all the analyses (R Core Team, 2017). Cluster analyses were run using the mclust package in R environment (Fraley, Raftery, Murphy, & Scrucca, 2012). For all the cluster analyses the default settings of mclust were used, more information about model comparisons, estimating procedures and starting values are provided in the technical publication (Scrucca, Fop, Murphy, & Raftery, 2016).

Graphs were obtained using the ggplot2 package in R environment (Wickham, 2009). Bootstraps were also performed using ggplot2, which is calculating 95% CI, based on the means, and using 1,000 boots. In all the analyses, both the statistical significance and the magnitude of the difference in terms of the effect size were presented.

**Results**

ANXIETY PROFILES AND PROTECTIVE FACTORS

Table 1 shows the correlations between all the measures and the children's age in months, together with the descriptive statistics. Interestingly, negative moderate correlations were found between the two SC scores (i.e., Academic and Competence subscales) and all anxiety measures (i.e., GA, TA and MA). Small negative correlations between ER scores and both GA and MA were observed. Finally, small negative correlations between the ABS scale and all anxiety measures (i.e., GA, TA and MA scales) were found.

Table 1 about here

We used a model-based clustering analysis approach, modelling clusters as a finite mixture of Gaussian distribution fitted via the EM algorithm (Fraley & Raftery, 2002). This method enabled us to assess different clustering solutions in terms of the model parameters. Several clustering models were considered, in terms of their evidence measured against the BIC index, and the number of underlying components (from 1 to 4). The model with 3 risk profiles proved to be superior (i.e., it had a higher BIC), and this solution was retained for subsequent analyses (see Table 2).

Table 2 about here

Model-based clustering of the GA, MA, and TA scores yielded three profiles (Figure 1). Profile 1 (N = 79, 11.72%) was labeled "low-risk" (because the anxiety scores were quite low overall); Profile 2 (N = 454, 66.22%) was labeled "average-risk" and refers to children with average levels of anxiety; finally, Profile 3 (N = 131, 22.06%) was labeled "high-risk" (because the levels of anxiety were high, particularly for GA and TA). Performances of GA, MA, and TA in these three profiles were not homogenous (Figure 1).

Figure 1 about here

A mixed ANOVA was performed, having anxiety as the dependent variable, with 3 risk profile [low, average and high] × 3 form of anxiety [GA, MA and TA] design was performed. The main effect of risk profile, $F(2, 661) = 391.77$, $p < .001$, $\eta^2_g = .351$, the main effect of form of anxiety, $F(2, 1322) = 113.18$, $p < .001$, $\eta^2_g = .085$, and the interaction effect, $F(4, 1322) = 33.93$, $p < .001$, $\eta^2_g = .0053$, were statistically significant[2]. All the post-hoc tests, with Bonferroni's corrections, were statistically significant ($p$s $< .05$) except for the difference between MA and GA or TA in the low-risk group ($p > .05$). The presence of a statistically significant interaction demonstrated that the profiles in the three groups were not flat. To further elucidate this finding we calculated the standardized difference between the three profiles across the three forms of anxiety (GA, MA, and TA). The results revealed extremely large effect sizes (*Cohen's d* $> 1.13$) in all cases except for the difference in MA between the average- and high-risk groups, which was small in terms of effect size (*Cohen's d* $= .46$), although statistically significant. Based on these findings, we surmise that the high-risk profile was characterized by higher levels of GA and TA, and lower levels of MA, while the average-risk profile featured average levels of GA and MA, and lower levels of TA.

**Differences between profiles on personal protective factors**

We ran a series of MANOVAs using school grade (grades 3 to 6) and risk profile (low, average and high risk) as fixed factors. In the first MANOVA, general protective factors (Self-Concept - Competence scale and Ego-Resiliency scale) were used as dependent variables. There was a statistically significant effect of the risk profile, $F(4, 1302) = 29.26$, $p < .001$; Wilk's $\Lambda = 0.842$, $\eta^2_p = .082$, but not of school grade, $F(6, 1302) = 0.891$, $p = .501$,

Wilk's $\Lambda = 0.992$, or the interaction between school grade and risk profile, $F(12, 1302) = 1.49$, $p = .121$, Wilk's $\Lambda = 0.973$, $\eta^2_p = .014$. Post-hoc ANOVAs confirmed a statistically significant effect of the risk profile on both Self-Concept – Competence, $F(2, 652) = 57.93$, $p < .001$, $\eta^2_p = .151$, and Ego-Resiliency, $F(2, 652) = 3.51$, $p = .031$, $\eta^2_p = .011$ (Figure 2).

In the second MANOVA, academic protective factors (Self-Concept – Academic scale and Academic Buoyancy Scale) were used as dependent variables. There was a statistically significant effect of the risk profile, $F(4, 1302) = 22.50$, p $< .001$; Wilk's $\Lambda = 0.875$, $\eta^2_p = .082$. The effect of school grade was statistically significant, $F(6, 1302) = 2.89$, $p = .008$, Wilk's $\Lambda = 0.974$, $\eta^2_p = .013$, while the interaction between school grade and risk profile, $F(12, 1302) = 0.71$, $p = .744$, Wilk's $\Lambda = 0.987$, $\eta^2_p = .006$, was not. Post-hoc ANOVAs confirmed a significant effect of the risk profile on both the Self-Concept – Academic scale, $F(2, 652) = 44.02$, $p < .001$, $\eta^2_p = .119$, and the Academic Buoyancy Scale, $F(2, 652) = 9.33$, $p < .001$, $\eta^2_p = .028$ (Figure 2). The effect of school grade was statistically significant for the Self-Concept – Academic scale, $F(2, 652) = 3.18$, $p < .024$, $\eta^2_p = .014$, but not for the Academic Buoyancy Scale, $F(2, 652) = 1.79$, $p = .147$, $\eta^2_p = .008$ (Figure 3).

A series of post-hoc analyses (using Bonferroni's correction) was also performed on the effects found statistically significant. In terms of the effect of risk profile, the three profiles differed from each other on the Self-Concept Competence ($p < .05$) and Academic scales ($p < .05$). The low- and average-risk profiles were statistically different on the Ego-Resiliency scale ($p < .05$), while the low- and high-risk profiles differed statistically on the Academic Buoyancy Scale ($p < .05$). As for the effect of school grade, children in grades 3 and 6 were statistically different on the Self-Concept – Academic scale ($p < .05$).

Figure 2 and 3 about here

## Discussion

The main aim of this study was to identify any presence of subgroups expressing particular anxiety patterns in a sample of school-age children, based on measures of GA, TA, and MA, using a latent profile analysis. We also tested the role of personal factors protecting against anxiety, distinguishing between general and academic self-concept and resilience.

Concerning our first aim, we found a three-profile solution: a small proportion of children (around 12%) in our sample expressed a low risk, meaning that they showed very low scores in various forms of anxiety; a large proportion (around 66%) showed an average risk, exhibiting average scores of different kinds of anxiety; and the remaining 22% revealed a high-risk profile of experiencing different form of anxiety. The profiles in these three groups were not flat, however: no differences emerged in the levels of GA, TA, and MA in the low-risk group, in which all children obtained low scores of anxiety on all measures; the average-risk profile was characterized by higher levels of GA and MA, and lower levels of TA, but in all cases the scores ranged around the average; and the high-risk profile featured higher levels of GA and TA, but lower levels of MA, although all scores ranged in high levels of anxiety.

Our latent profile analysis revealed different patterns of anxiety, albeit without any clear distinction between general and academic anxiety (which only becomes apparent in older students; on this point, see also Campbell & Rapee, 1996). According to Carey and colleagues (2017), MA can develop through two main mechanisms: a predisposition to anxiety in general, or a repeatedly poor performance in mathematics. Our data seem to support the hypothesis that MA in younger children could be driven primarily by a general tendency to be anxious, since our average-risk profile was characterized by higher levels of GA and MA than of TA. Our findings do not exactly replicate those of Carey and colleagues (2017), however, partly because they tested students in grades 4 or 7-8, while our children were in grades 3 to 6. In fact, the effect of age was not statically significant in our sample, so we found neither a simple distinction between different (i.e., low, average, high) degrees of risk to develop anxiety, nor any clear difference between general and academic forms of anxiety. In particular, testing children along a continuum, from grades 3 to 6, enabled us to identify an intermediate stage of anxiety development, at which point our high-risk profile coincided with more GA and TA than MA, while our average-risk profile showed more GA and MA than TA. Judging from our data, the hypothesis advanced by Carey and colleagues (2017) for the development of MA may therefore extend to TA as well. In other words, GA seems to act as a risk factor for the onset of other, more or less severe forms of anxiety (MA and TA).

To better analyze developmental changes, further studies should try to replicate our findings in older students too. In particular, longitudinal research could help to   clarify whether the forms of anxiety seen in younger children are precursors of those identified in older students. The development of different forms of anxiety has so far been underexplored in the field of psychological research, and more information on this topic would help us to develop new, effective interventions, adapting their content to a given child's developmental age and specific anxiety profile.

The second goal of our study was to test the effects of personal general and academic factors protecting against anxiety. We found that distinguishing between general and academic protective factors did not produce different results by risk profile among children: the effect of self-confidence, be it general or academic, differentiated between the three profiles more clearly  than the effect of resilience or academic buoyancy. In particular, children with a low-risk profile reported feeling more competent and had a stronger academic self-concept than children with average- or high-risk profiles. As for age differences, the correlations between age and personal protective factors were small. An intriguing effect emerged on the Self-Concept – Academic scale, however, which was significantly higher for 3[rd]-graders than for 6[th]-graders. The existing literature offers ample evidence of high levels of academic self-concept predicting low levels of TA (Bong & Skaalvik, 2002; Putwain & Daniels, 2010), and shown that MA promotes negative academic self-concepts regarding math abilities (Ashcraft & Kirk, 2001; Wu et al., 2012). Bandalos and colleagues (1995) also found that academic self-concept mediated the relation between prior experience and anxiety. Our study is the first to have jointly considered the relationship between different forms of anxiety and both academic and general self-concept. Self-competence was clearly inversely proportional to the level of anxiety: the higher the risk of becoming anxious, the lower the self-concept (on both the competence and the academic scales). Academic self-concept also tended to decline as the children grew older age, giving the impression that an older age (and possibly more negative academic experiences), is associated to a higher risk of anxiety, and to a worsening specific academic self-concept.

As for the effect of resilience, we found that school children with a low anxiety risk profile scored higher than those at higher risk on both measures of resilience and academic buoyancy. The children at average risk of anxiety did not differ from those with a high-risk profile (or from the low-risk group when it came to academic buoyancy). Consistently with the definition of this construct, school children are therefore able to maintain a certain level of resilience at least up until they develop an average risk of becoming anxious. In fact, resilience only decreased in our high-risk group (with higher levels of both GA and TA). Thus, resilience and academic buoyancy are both personal factors crucial in  protecting against anxiety. Resilient students can maintain or regain prior levels of functioning, unlike those students who respond badly to adversity, and consequently risk developing higher levels of anxiety.

Although it contains some very interesting findings this study also has some limitations that should be addressed in future publications. First, our agreement with the schools did not include testing neither of first and second graders nor of students older than grade 8, hence further studies should also consider younger children in order to better understand the development of different forms of anxiety. Moreover, future studies should evaluate the role of general and academic personal protective factors in older children. This would be very important for understanding the complex interlink between these variables during the development. Second, our instruments included several self-reported measures, which can be somewhat affected by common-method bias (Williams & Brown, 1994), and were not administered in a counterbalanced order, similarly as already done in previous studies (e.g., Giofrè, Borella & Mammarella, 2017; Justicia-Galiano, et al., 2017; Hill et al., 2016). Further research should replicate our outcomes also controlling for these small methodological issues.

It is worth noting that our findings also have important clinical and educational implications: for a start, a better understanding of how different forms of anxiety develop is crucial to design specific intervention programs aimed at reducing students' anxiety; as well as, fostering resilience can prevent academic anxiety. Our findings can be used to underpin the preparation and testing of interventions (Rose, Miller, & Martinez, 2009). In particular, interventions could be focused in promoting academic self-concept and resilience offering continuous feedback, suitable tasks based on children's different competence levels, and promoting collaborative learning (Hagenauer & Hascher 2014). In addition, school-based interventions, based on cognitive and behavioral emotion regulation strategies for specific forms of anxiety, could reduce the onset of anxiety symptoms (Weems et al. 2015) and at the same time foster protective factors: students with high levels of self-concept and resilience are able to take on challenging tasks, persist when faced with difficulties, and believe in their ability to do well (Niemiec & Ryan, 2009).

To conclude, our study suggest that forms of anxiety start to differentiate between 3[rd] and 6[th] grade. Although no clear distinction between general and academic anxiety emerged in this age bracket, GA seemed to act as a risk factor for the onset of more specific forms of anxiety. Concerning the role of personal protective factors, different effects emerged for self-concept and resilience (both general and academic): self-concept was lower the higher the anxiety risk profile (from low to average and high); and academic self-concept also decreased with age, while resilience only decreased in association with the high anxiety risk profile.

Press.

**Footnotes**

[1] A demographic questionnaire was administered in order to estimate the socio-economic condition of the family. In particular, children were required to report whether their parents were employed or not and to briefly describe the job of both parents.

[2] The analysis was also repeated including school grades. The effect of grade, the two-way interactions between grade and risk profile, and between grade and anxiety, and the three-way interaction between grade, risk profile and anxiety were small and not statistically significant ($ps > .051$, $\eta^2_g < .004$).

**Acknowledgement**

The present study was supported by BIRD161109/16, funded by the University of Padova (Italy).

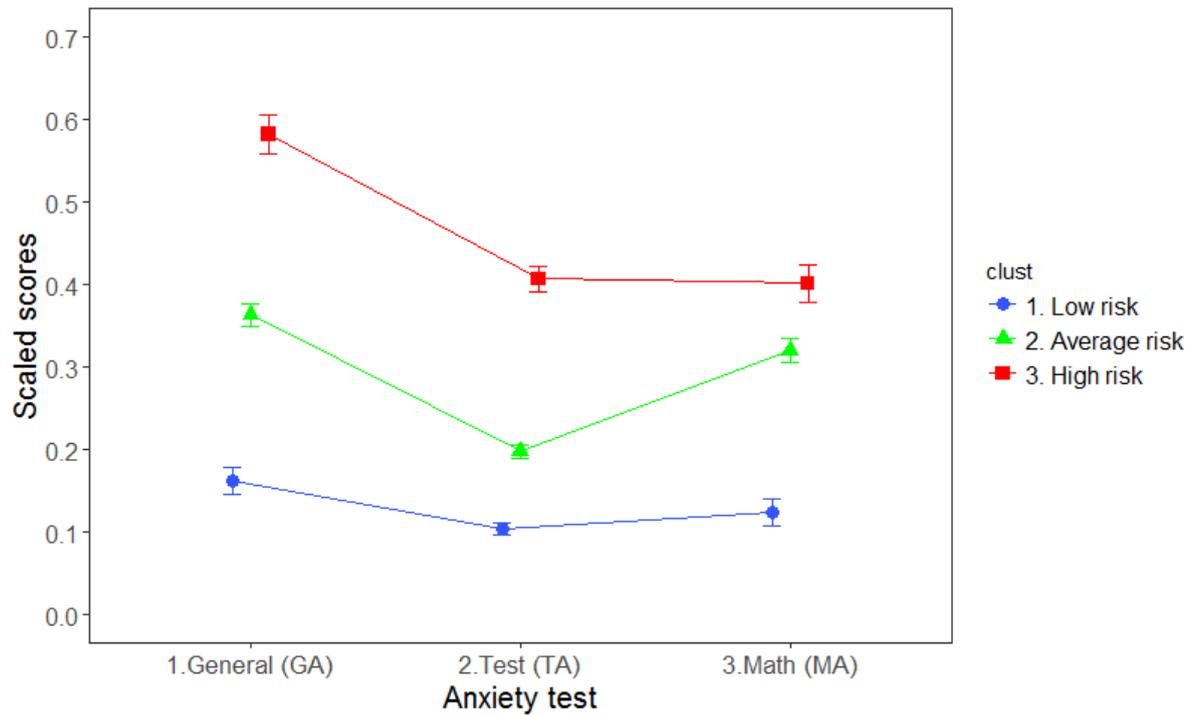

*Figure 1*. The profile of the three clusters - low-, average- and high-risk - on general anxiety (GA), test anxiety (TA) and math anxiety (MA). Error bars represent 95% bootstrapped confidence intervals of the interaction.

s

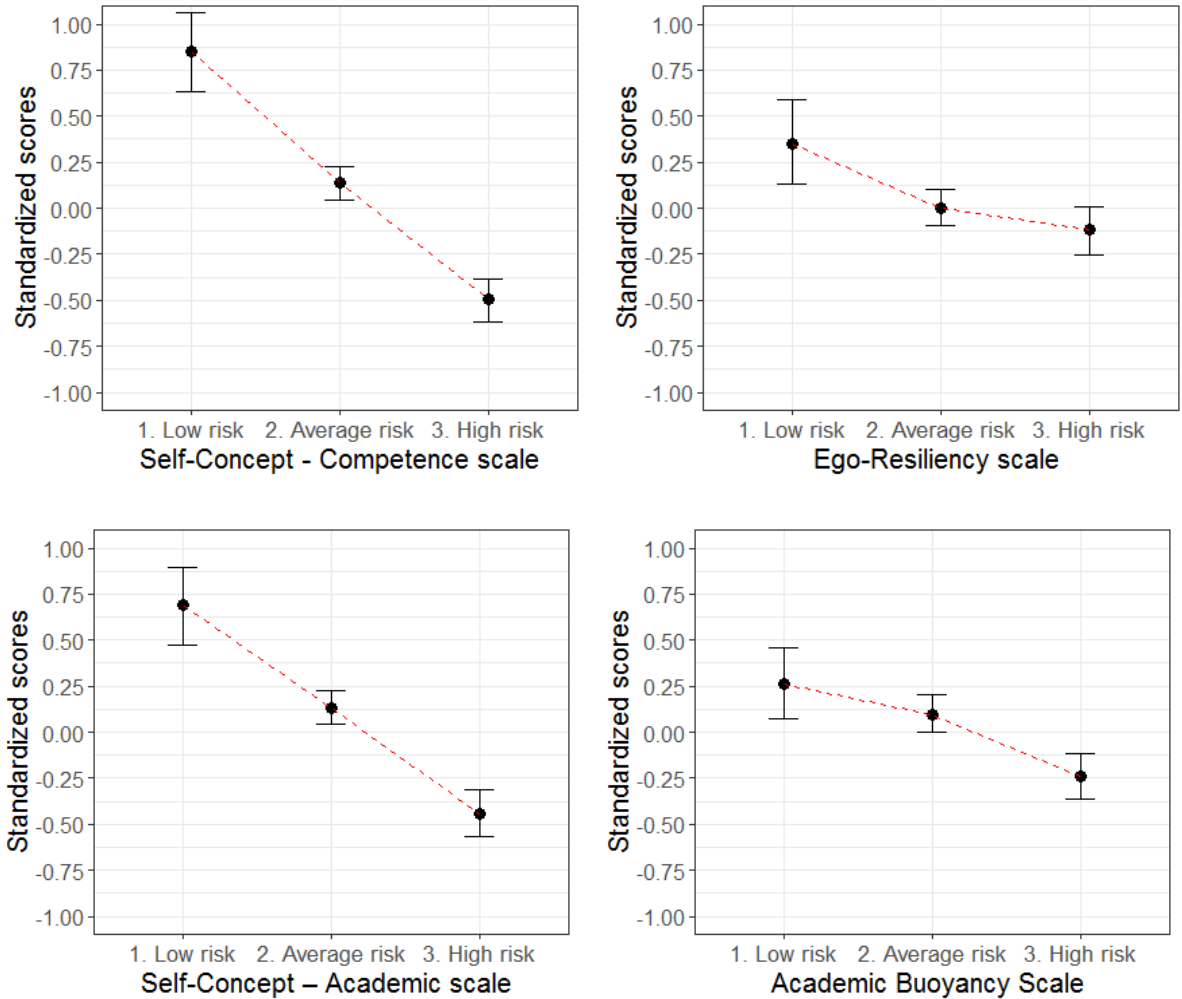

*Figure 2*. Performance of the three anxiety profiles (low-, average- and high-risk) on general (above) and academic (below) protective factors. Higher scores represent better self-concept, higher resilience and higher academic buoyancy. Standard errors represent bootstrapped 95% Confidence Intervals.

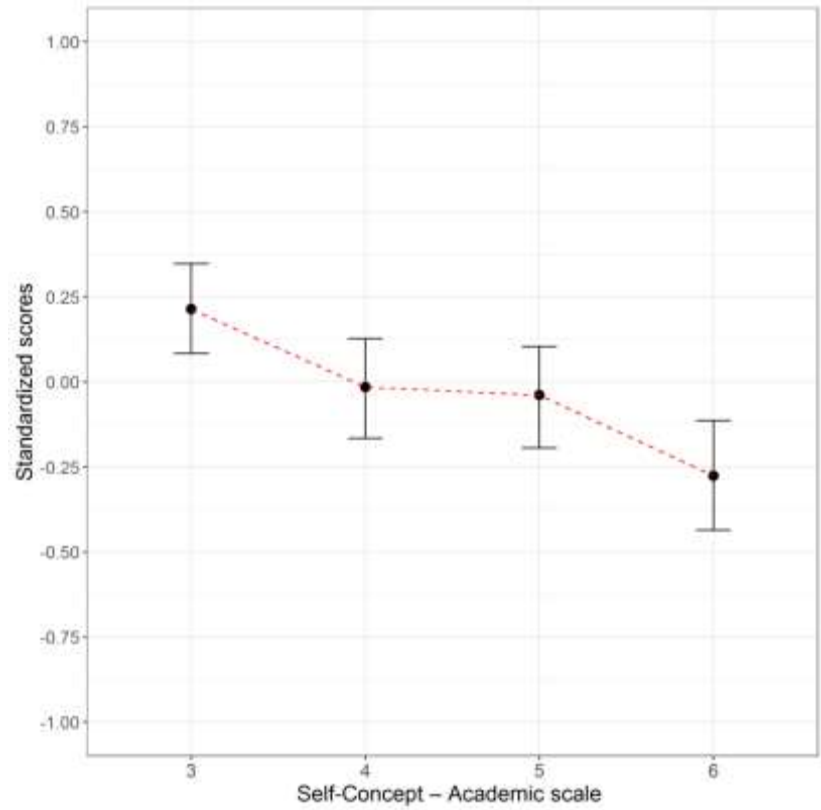

*Figure 3*. Effects of school grade (3ʳᵈ, 4ᵗʰ, 5ᵗʰ and 6ᵗʰ) on Self-concept - Academic scale. Error

bars represent bootstrapped confidence intervals.

Table 1

*Correlations, means and standard deviations for the sample as a whole.*

|  | 1 | 2 | 3 | 4 | 5 | 6 | 7 | 8 |
|---|---|---|---|---|---|---|---|---|
| 1. Age | 1 |  |  |  |  |  |  |  |
| 2. General anxiety | .031 | 1 |  |  |  |  |  |  |
| 3. Test anxiety | .068 | .596[*] | 1 |  |  |  |  |  |
| 4. Math anxiety | .056 | .413[*] | .337[*] | 1 |  |  |  |  |
| 5. SC-Competence | -.083[*] | -.444[*] | -.427[*] | -.282[*] | 1 |  |  |  |
| 6. ER | .018 | -.177[*] | -.075 | -.123[*] | .354[*] | 1 |  |  |
| 7. SC-Academic | -.152[*] | -.418[*] | -.404[*] | -.325[*] | .668[*] | .365[*] | 1 |  |
| 8. ABS | .026 | -.122[*] | -.151[*] | -.157[*] | .248[*] | .237[*] | .250[*] | 1 |
| M | 9.20 | 15.39 | 54.46 | 22.66 | 74.54 | 4.62 | 72.29 | 19.83 |
| SD | 1.13 | 7.58 | 14.90 | 7.80 | 1.00 | 6.22 | 9.38 | 5.87 |

*Note*. SC-Competence = self-concept competence scale; ER= Ego-resiliency; SC-Academic =

self-concept academic scale; ABS = Academic buoyancy scale

[*]$p < .05$

Table 2

*BIC for clustering models as function of the number of components. The higher the BIC values, the better the model*

| Number of components | BIC values | |
| :---: | :---: | :---: |
| | VEE | VVE |
| 1 | 1754.48 | 1754.48 |
| 2 | 1812.02 | 1807.15 |
| **3** | **1822.06** | **1812.57** |
| 4 | 1787.66 | 1761.06 |
| 5 | 1769.45 | 1743.65 |
| 6 | 1764.25 | 1724.45 |
| 7 | 1766.50 | 1701.52 |
| 8 | 1736.84 | 1668.23 |
| 9 | 1721.15 | 1641.95 |

*Note.* VEE = ellipsoidal, equal shape and orientation; VVE = ellipsoidal, equal orientation.

Higher BIC values correspond to a better fit. For further information, see Scrucca et al. (2016),

the maximum BIC value occurs for the three clusters model, which is in bold in the table.